\documentclass[fleqn,10pt]{wlscirep}
\usepackage{algorithm}
\usepackage{algorithmic}
\usepackage{bm}
\usepackage[T1]{fontenc}
\usepackage{here}
\usepackage[utf8]{inputenc}
\usepackage{mathtools}
\usepackage{subcaption}
\usepackage{comment}

\title{Estimating profitable price bounds for prescriptive price optimization}

\author[1]{Masato Inokuma}
\author[2,*]{Shunnosuke Ikeda}
\author[2,3]{Yuichi Takano}
\affil[1]{Graduate School of Science and Technology, University of Tsukuba, Tsukuba, Japan}
\affil[2]{Institute of Systems and Information Engineering, University of Tsukuba, Tsukuba, Japan}
\affil[3]{Center for Artificial Intelligence Research, Tsukuba Institute for Advanced Research (TIAR), University of Tsukuba, Tsukuba, Japan}

\affil[*]{ikeda@cs.tsukuba.ac.jp}


\begin{abstract}
Pricing of products and services, which has a significant impact on consumer demand, is one of the most important factors in maximizing business profits.
Prescriptive price optimization is a prominent data-driven pricing methodology consisting of two phases: demand forecasting and price optimization. 
In the practice of prescriptive price optimization, the price of each item is typically set within a predetermined range defined by lower and upper bounds. 
Narrow price ranges can lead to missed opportunities, while wide price ranges run the risk of proposing unrealistic prices; therefore, determining profitable price bounds while maintaining the reliability of the suggested prices is a critical challenge that directly affects the effectiveness of prescriptive price optimization.
We propose two methods for estimating price bounds in prescriptive price optimization so that future total revenue derived from the optimized prices will be maximized. 
Our first method for price bounds estimation uses the bootstrap procedure to estimate confidence intervals for optimal prices. 
Our second method uses the Nelder–Mead simplex method for black-box price bounds optimization that maximizes total revenue estimated through $K$-fold cross-validation.
Experimental results with synthetic price--demand datasets demonstrate that our methods successfully narrowed down the price range while maintaining high revenues, particularly when the number of items was small or the demand noise level was low. 
Moreover, as more data accumulated, the comparative advantage of our methods further increased.
\end{abstract}
\begin{document}

\flushbottom
\maketitle
%
%
\thispagestyle{empty}


\section*{Introduction}
\subsection*{Background}

Pricing of products and services, which has a significant impact on consumer demand, is one of the most important factors in maximizing business profits~\cite{bu2023offline}.
In fact, appropriate pricing strategies can bring large benefits to a variety of businesses~\cite{shams2020price}, such as every day low prices in retail stores~\cite{mehrotra2020price}, price promotions for beauty category products~\cite{phumchusri2024price}, and strategic revenue management in restaurants~\cite{tyagi2022approaches}. 
In addition, recent advances in information and communication technology have made it possible to collect huge amounts of consumer purchase data, enabling the development of effective data-driven pricing strategies.
As a result, the impact of pricing strategies on business revenues and profits is more pronounced than ever before~\cite{shams2020price}.

Against this background, \emph{prescriptive price optimization} \cite{ito2017optimization,ikeda2023prescriptive} has emerged as a prominent data-driven pricing methodology.
This framework consists of two sequential phases: demand forecasting and price optimization.
In the first demand forecasting phase, regression models are built to quantify the impact of prices on consumer demand for multiple items (e.g., products or services).
In the subsequent price optimization phase, these demand forecasting models are used to formulate and solve a price optimization problem for maximizing a specified objective function (e.g., total revenue or profit).

In the practice of prescriptive price optimization, the price of each item is typically set within a predetermined range defined by lower and upper bounds.
Setting these ranges too narrowly may prevent large price changes and lead to missed opportunities.
Conversely, setting price ranges too broadly allows for large price changes but at the same time runs the risk of proposing unrealistic prices.
Accordingly, determining profitable price bounds while maintaining the reliability of the suggested prices is a critical challenge that directly affects the effectiveness of prescriptive price optimization.

\subsection*{Related work}
Several studies have been conducted on prescriptive price optimization from various perspectives~\cite{ito2017optimization, ikeda2023prescriptive, ferreira2016analytics, ito2016large, yabe2017robust, ito2018unbiased}.
The major challenges in prescriptive price optimization include: 
(i) high computational complexity depending on the demand forecasting model and optimization algorithm employed, and
(ii) overestimation of revenues and profits caused by errors in demand forecasts---a phenomenon known as \emph{optimistic bias}~\cite{smith2006optimizer}.
The former challenge makes it difficult to calculate the optimal prices, while the latter challenge imperils the validity of the calculated optimal prices.

Various modeling and solution methods have been proposed to reduce the computational complexity of prescriptive price optimization. 
Ferreira et al.~\cite{ferreira2016analytics} employed an algorithm based on linear relaxation of the price optimization problem using bagging-based regression trees for demand forecasting. 
Ito and Fujimaki~\cite{ito2016large} designed a network flow algorithm for large-scale price optimization based on the relationship between revenue supermodularity and demand cross-elasticity. 
Ito and Fujimaki~\cite{ito2017optimization} developed a fast approximation algorithm using semidefinite relaxation of the price optimization problem.
Ikeda et al.~\cite{ikeda2023prescriptive} formulated the price optimization problem using optimal regression trees for demand forecasting and developed a heuristic algorithm based on the randomized coordinate ascent. 

To mitigate overestimation of revenues and profits in prescriptive price optimization, Yabe et al.~\cite{yabe2017robust} devised a robust quadratic optimization model that accounts for the uncertainty of demand forecasts represented by matrix normal distribution. 
Ito et al.~\cite{ito2018unbiased} proposed two methods for unbiased objective estimation: the cross-validation-based method and the parameter perturbation method. 
They reported that these two methods successfully corrected the optimistic bias on both synthetic and real-world datasets.

In these prior studies, optimal prices have been determined subject to the constraints of lower and upper price bounds.
Similar constraints have been imposed in most studies on price optimization (e.g., robust price optimization \cite{thiele2006single, hamzeei2022robust} and dynamic pricing \cite{harsha2019practical, shao2025constrained}); see the comprehensive reviews \cite{bitran2003overview, soon2011review} for detailed discussions on multi-item price optimization models. 
In practical applications, appropriate price bounds vary depending on product type and sales period, and also have a significant impact on revenues and profits.
None of the aforementioned studies, however, have specifically examined how to determine appropriate price bounds for prescriptive price optimization.

Recently, Ikeda et al.~\cite{ikeda2023operating,ikeda2024interpretable} proposed a comprehensive framework for setting price bounds in prescriptive price optimization. 
This framework consists of estimating and adjusting price bounds; specifically, price bounds are first estimated from historical pricing data, and then these estimated price bounds are adjusted by solving a shape-constrained optimization problem.
However, this framework is focused on the interpretability of optimal prices and their consistency with observed data, and does not aim to maximize total revenue or profit generated by price optimization.


\subsection*{Our contribution}

The goal of this paper is to establish a computational framework for estimating price bounds for revenue maximization in prescriptive price optimization. 
First, we present a problem formulation for price bounds estimation by following bilevel optimization approaches to machine learning model selection~\cite{bennett2008bilevel,sinha2017review,takano2020best}. 
Next, we propose two methods for estimating profitable
price bounds; the first method uses the bootstrap procedure~\cite{efron1992bootstrap,james2013introduction} to estimate confidence intervals for optimal prices, and the second method uses the Nelder–Mead simplex
method~\cite{nelder1965simplex,kochenderfer2019algorithms} for black-box price bounds optimization to maximize total revenue estimated through $K$-fold cross-validation~\cite{ito2018unbiased}.

To demonstrate the effectiveness of our two methods for price bounds estimation, we conducted numerical experiments using synthetic price--demand datasets.
Our methods successfully narrowed down the price range while maintaining high revenues, particularly when the number of items was small or the demand noise level was low.
As more data accumulated, the comparative advantage of our methods further increased.
Additionally, the price ranges calculated by our two methods showed contrasting trends depending on the accuracy of demand forecasts.
The computation time for the bootstrap estimation method increased only slightly with the number of items, whereas that for the cross-validation-guided black-box optimization method was significantly affected by the item count.

\section*{Methods}

In this section, we first give an overview of the framework of prescriptive price optimization~\cite{ito2017optimization,ikeda2023prescriptive}. 
We next present our ideal problem formulation to clarify the target problem of estimating price bounds for revenue maximization. 
We then propose two methods for estimating profitable price bounds: bootstrap estimation and cross-validation-guided black-box optimization. 

\subsection*{Prescriptive price optimization}
The framework of prescriptive price optimization~\cite{ito2017optimization,ikeda2023prescriptive} consists of two sequential phases: demand forecasting and price optimization. 
For simplicity, this paper is aimed at maximizing total revenue. 
We also denote the set of first $m$ positive integers as $[m] \coloneqq \{1,2,\dots,m\}$.

Let $\bm{p} \coloneqq (p_j)_{j \in [m]} \in \mathbb{R}^{m}$ denote a vector composed of $m$ item prices.
The demand for each item $j \in [m]$ is expressed as a function of the item prices, $d_j (\bm{p},\bm{\theta}_j^\star)$, where $\bm{\theta}_j^\star$ is a vector of function parameters. 
Total revenue is given by the sum of the products of demand and price as 
\begin{align}\label{eq:total_revenue}
    f(\bm{p},\bm\Theta^\star) \coloneqq \sum_{j = 1}^m d_j (\bm{p},\bm{\theta}_j^\star) p_j,
\end{align}
where $\bm\Theta^\star \coloneqq (\bm{\theta}_j^\star)_{j \in [m]}$ is a parameter matrix.
The objective of prescriptive price optimization is to find a price vector $\bm p \in \mathbb{R}^{m}$ that maximizes the total revenue (Eq.~\eqref{eq:total_revenue}).

In practice, however, the ground-truth parameter matrix $\bm\Theta^\star$ is unknown.
Accordingly, we instead use a parameter matrix $\hat{\bm{\Theta}}$ estimated from historical data.
Assuming a linear demand function of prices, two sequential phases of prescriptive price optimization are described as follows.

\paragraph{(i) Demand forecasting:}
Estimate demand forecasting functions for each item $j \in [m]$ as
\begin{align}\label{eq:demand_function}
    d_j (\bm{p},\hat{\bm{\theta}}_j) \coloneqq \hat\theta_{j0} + \sum_{\ell=1}^m \hat\theta_{j \ell} p_\ell \quad (j \in [m]),
\end{align}
where $\hat{\bm{\theta}}_j \coloneqq (\hat{\theta}_{j\ell})_{\ell \in \{0\}\cup[m]} \in \mathbb{R}^{m+1}$ is a vector of regression coefficients for each item $j \in [m]$.
The corresponding coefficient matrix $\hat{\bm{\Theta}} \coloneqq (\hat{\bm{\theta}}_j)_{j \in [m]} \in \mathbb{R}^{(m+1) \times m}$ is obtained through ordinary least squares estimation using a historical price--demand dataset $\bm{X} \coloneqq \{(\bm{p}_i,\bm{d}_i) \mid i \in [n]\}$, where $\bm{p}_i \coloneqq (p_{ij})_{j \in [m]} \in \mathbb{R}^{m}$ and $\bm{d}_i \coloneqq (d_{ij})_{j \in [m]} \in \mathbb{R}^{m}$ denote the price and demand vectors for each data instance $i \in [n]$, respectively. 

\paragraph{(ii) Price optimization:}
Solve the following price optimization problem based on the estimated coefficient matrix $\hat{\bm{\Theta}} \in \mathbb{R}^{(m+1) \times m}$ for revenue maximization:
\begin{align}
    \underset{\bm{p} \in \mathbb{R}^m}{\text{maximize}}\quad& f(\bm{p},\hat{\bm{\Theta}}) \coloneqq \sum_{j = 1}^m d_j (\bm{p},\hat{\bm{\theta}}_j) p_j \label{eq:obj}\\
    \text{subject to}\quad&\bm{\alpha} \leq \bm{p} \leq \bm{\beta}, \label{eq:constr_price}
\end{align}
where $\bm{\alpha} \coloneqq (\alpha_j)_{j \in [m]} \in \mathbb{R}^m$ and $\bm{\beta} \coloneqq (\beta_j)_{j \in [m]} \in \mathbb{R}^m$ represent the lower and upper price bounds, respectively.


\subsection*{Ideal bilevel formulation for price bounds estimation}

Note that solving the price optimization problem (Eq.~\eqref{eq:obj}--\eqref{eq:constr_price}) does not necessarily lead to maximizing the ground-truth total revenue (Eq.~\eqref{eq:total_revenue}) due to estimation errors in the parameter matrix $\hat{\bm{\Theta}} \in \mathbb{R}^{(m+1) \times m}$. 
As such, we focus on estimating price bounds such that solving the problem (Eq.~\eqref{eq:obj}--\eqref{eq:constr_price}) improves the ground-truth total revenue (Eq.~\eqref{eq:total_revenue}). 
 
Let $\bm{p}^{\mathrm{min}} \coloneqq (p^{\mathrm{min}}_j)_{j \in [m]} \in \mathbb{R}^{m}$ and $\bm{p}^{\mathrm{max}} \coloneqq (p^{\mathrm{max}}_j)_{j \in [m]} \in \mathbb{R}^m$ denote the vectors of minimum and maximum feasible prices, respectively.
Estimating the price bounds for maximizing the ground-truth total revenue (Eq.~\eqref{eq:total_revenue}) can be posed as the following bilevel optimization problem:
\begin{align}
    \underset{\hat{\bm{p}},\bm{\alpha},\bm{\beta} \in \mathbb{R}^m}{\text{maximize}}\quad& f(\hat{\bm{p}},\bm{\Theta}^{\star}) \coloneqq \sum_{j = 1}^m d_j (\hat{\bm{p}},\bm{\theta}^\star_j) \hat{p}_j \label{eq:obj2}\\
    \text{subject to}\quad&\hat{\bm p}
      \;\in\;
      \underset{\bm p \in \mathbb{R}^m} {\operatorname{argmax}}
      \Bigl\{ f(\bm p,\hat{\bm\Theta})
        \;\Bigm|\;
        \bm{\alpha} \le \bm p \le \bm{\beta}
      \Bigr\}, \label{eq:lower}\\
      & \sum^m_{j = 1}(\beta_j - \alpha_j) \leq \Gamma, \label{eq:width}\\
      & \bm p^{\min} \le \bm{\alpha} \le \bm{\beta} \le  \bm p^{\max}, \label{eq:bounds}
\end{align}
where $\Gamma \in \mathbb{R}_{+}$ is a parameter for controlling the total width of the price range as in Eq.~\eqref{eq:width}.
Note here that $\hat{\bm{p}} \coloneqq (\hat{p}_j)_{j \in [m]} \in \mathbb{R}^{m}$ is a solution to the price optimization problem (Eq.~\eqref{eq:obj}--\eqref{eq:constr_price}) with the price bounds $\bm{\alpha},\bm{\beta} \in \mathbb{R}^m$.
However, it is impossible to directly solve the bilevel optimization problem (Eq.\eqref{eq:obj2}--\eqref{eq:bounds}) because the ground-truth coefficient matrix $\bm\Theta^\star \in \mathbb{R}^{(m+1) \times m}$ in Eq.~\eqref{eq:obj2} is unknown in practice. 

\subsection*{Bootstrap estimation}
Our first method for estimating profitable price bounds involves using the bootstrap procedure \cite{efron1992bootstrap,james2013introduction} to estimate confidence intervals for the optimal prices as an estimate of the lower and upper price bounds.
Algorithm~\ref{alg:bs} describes our bootstrap method for estimating profitable price bounds. 

This algorithm first generates bootstrap samples $\bm{X}^{(b)}$ for $b \in [N^{{\textrm bs}}]$ by repeating $N^{{\textrm bs}}$ times the process of extracting $n$ instances with replacement from the historical dataset $\bm{X} = \{(\bm{p}_i,\bm{d}_i) \mid i \in [n]\}$.
From each bootstrap sample $\bm{X}^{(b)}$, we estimate the coefficient matrix $\hat{\bm{\Theta}}^{(b)} \in \mathbb{R}^{(m+1) \times m}$ of demand forecasting functions (Eq.~\eqref{eq:demand_function}) through the ordinary least squares estimation. 
We then compute the corresponding optimal prices $\hat{\bm{p}}^{(b)} \coloneqq (\hat{p}^{(b)}_j)_{j \in [m]} \in \mathbb{R}^m$ by solving problem (Eq.~\eqref{eq:obj}--\eqref{eq:constr_price}) based on the estimated coefficient matrix $\hat{\bm{\Theta}}^{(b)} \in \mathbb{R}^{(m+1) \times m}$.
For each item $j\in[m]$, we calculate the sample mean $\bar{p}_j$ and sample standard deviation $s_{j}$ as in Eq.~\eqref{eq:se} from the obtained collection of optimal prices. 
Finally, we set symmetric price bounds as in Eq.~\eqref{eq:kappa}, where $\kappa \in \mathbb{R}_{+}$ is a critical-value parameter for controlling the width of each price range. 

\begin{algorithm}[thb!]
\caption{Bootstrap method for estimating profitable price bounds}
\label{alg:bs}
\begin{algorithmic}[1]
\item[\textbf{Input:}] 
Historical dataset $\bm{X} = \{(\bm{p}_i,\bm{d}_i) \mid i \in [n]\}$, 
number of bootstrap samples $N^{\mathrm{bs}} \in \mathbb{N}$, 
critical-value parameter $\kappa \in \mathbb{R}_{+}$, and
minimum and maximum feasible prices $\bm{p}^{\mathrm{min}}, \bm{p}^{\mathrm{max}} \in \mathbb{R}^m$.
    \FOR{$b \in [N^{{\textrm bs}}]$}
    \STATE Draw a bootstrap sample $\bm{X}^{(b)}$ from the dataset $\bm{X}$.
    \STATE Estimate the coefficient matrix $\hat{\bm{\Theta}}^{(b)} \in \mathbb{R}^{(m+1) \times m}$ from the bootstrap sample $\bm{X}^{(b)}$. 
    \STATE Compute the optimal prices $\hat{\bm{p}}^{(b)} \in \mathbb{R}^m$ to problem (Eq.~\eqref{eq:obj}--\eqref{eq:constr_price}) based on the coefficient matrix $\hat{\bm{\Theta}}^{(b)} \in \mathbb{R}^{(m+1) \times m}$.
    \ENDFOR
    \STATE Calculate the sample mean and sample standard deviation of the optimal prices as 
    \begin{align}\label{eq:se}
        &\bar{p}_{j} \coloneqq \frac{1}{N^{{\rm bs}}}\sum_{b=1}^{N^{{\rm bs}}} \hat{p}^{(b)}_j, \quad s_{j} \coloneqq \sqrt{\frac{1}{N^{{\rm bs}}-1}\sum_{b=1}^{N^{{\rm bs}}} \biggl(\hat{p}^{(b)}_j - \bar{p}_{j})\biggr)^2} \quad (j \in [m]).
    \end{align}
    \STATE Calculate the lower and upper price bounds as
    \begin{align}\label{eq:kappa}
        \alpha_{j} \coloneqq \max\{p^{\mathrm{min}}_j, ~\bar{p}_{j} - \kappa s_{j}\}, \quad \beta_{j} \coloneqq \min \{p^{\mathrm{max}}_j, ~\bar{p}_{j} + \kappa s_{j} \}\quad (j \in [m]).
    \end{align}
\item[\textbf{Output:}] Lower and upper price bounds $\bm{\alpha},\bm{\beta} \in \mathbb{R}^m$. 
\end{algorithmic}
\end{algorithm}

\subsection*{Cross-validation-guided black-box optimization}
Our second method for estimating profitable price bounds involves using the $K$-fold cross-validation procedure proposed by Ito et al.~\cite{ito2018unbiased} for unbiased revenue estimation in prescriptive price optimization.
We then apply the Nelder--Mead simplex method~\cite{nelder1965simplex,kochenderfer2019algorithms} to search for the price bounds that maximize the estimated total revenue.

Algorithm~\ref{alg:kfoldCV} describes the detailed procedure of the $K$-fold cross-validation for estimating the ground-truth total revenue. 
This algorithm first divides the dataset $\bm{X} = \{(\bm{p}_i,\bm{d}_i) \mid i \in [n]\}$ into $K$ disjoint sub-datasets $\bm{X}_k$ for $k \in [K]$.
For each fold $k \in [K]$, we estimate the coefficient matrix $\bm\Theta_k^{\mathrm{trn}} \in \mathbb{R}^{(m+1) \times m}$ from the training dataset $\bm{X}\setminus\bm{X}_k$ and the coefficient matrix $\bm\Theta_k^{\mathrm{vld}} \in \mathbb{R}^{(m+1) \times m}$ from the validation dataset $\bm{X}_k$ through ordinary least squares estimation.
We compute the price vector $\hat{\bm{p}}_k^{\mathrm{trn}} \in \mathbb{R}^m$ using the training coefficient matrix $\bm{\Theta}_k^{\mathrm{trn}} \in \mathbb{R}^{(m+1) \times m}$ as in Eq.~\eqref{prob1trn}, and then evaluate the resultant total revenue using the validation coefficient matrix $\bm{\Theta}_k^{\mathrm{vld}} \in \mathbb{R}^{(m+1) \times m}$. 
By averaging the validation revenues over $k \in [K]$, we obtain $\mathrm{CV}_{\bm X}(\bm \alpha,\bm \beta)$ as an estimate of the ground-truth total revenue as in Eq.~\eqref{eq:cv}.

\begin{algorithm}[thb!]
\caption{$K$-fold cross-validation for estimating the ground-truth total revenue}
\label{alg:kfoldCV}
\begin{algorithmic}[1]
  \renewcommand{\algorithmicrequire}{\textbf{Input:}}
  \renewcommand{\algorithmicensure}{\textbf{Output:}}
  \REQUIRE 
  Historical dataset $\bm{X} = \{(\bm{p}_i,\bm{d}_i) \mid i \in [n]\}$, 
  number of folds $K$, and 
  lower and upper price bounds $\bm{\alpha},\bm{\beta} \in \mathbb{R}^m$. 
  \STATE Divide the dataset $\bm X$ into $K$ disjoint subsets $\bm X_k$ for $k \in [K]$.
  \FOR{$k \in [K]$}
    \STATE Estimate the coefficient matrix $\bm\Theta_k^{\mathrm{trn}} \in \mathbb{R}^{(m+1) \times m}$ from the training dataset $\bm X \setminus \bm X_k$.
    \STATE Estimate the coefficient matrix $\bm\Theta_k^{\mathrm{vld}} \in \mathbb{R}^{(m+1) \times m}$ from the validation dataset $\bm X_k$.
    \STATE Compute the optimal prices based on the training coefficient matrix $\bm\Theta_k^{\mathrm{trn}} \in \mathbb{R}^{(m+1) \times m}$ as follows:
    \begin{align} \label{prob1trn}
      \hat{\bm p}_k^{\mathrm{trn}}
      \;\in\;
      \underset{\bm p \in \mathbb{R}^m} {\operatorname{argmax}}
      \Bigl\{ f(\bm p,\bm\Theta_k^{\mathrm{trn}})
        \;\Bigm|\;
        \bm \alpha \le \bm p \le \bm \beta
      \Bigr\}. 
    \end{align}
  \ENDFOR
  \ENSURE The average of $K$ revenue estimates based on the validation coefficient matrix $\bm\Theta_k^{\mathrm{vld}} \in \mathbb{R}^{(m+1) \times m}$ for $k \in [K]$:
  \begin{align}
    \mathrm{CV}_{\bm X}(\bm \alpha,\bm \beta)
    \coloneqq \frac{1}{K} \sum_{k=1}^K
      f\bigl(\hat{\bm p}_k^{\mathrm{trn}},\bm\Theta_k^{\mathrm{vld}}\bigr).  \label{eq:cv}    
  \end{align}
\end{algorithmic}
\end{algorithm}

We aim to calculate the price bounds that maximize the estimated total revenue through cross-validation as 
\begin{align}
    \underset{\bm{\alpha},\bm{\beta} \in \mathbb{R}^m}{\text{maximize}} \quad &\mathrm{CV}_{\bm{X}}(\bm{\alpha},\bm{\beta}) \label{eq:obj_prob1}\\
    \text{subject to} \quad &\sum^m_{j = 1}(\beta_j - \alpha_j) \leq \Gamma, \label{con1:prob1}\\
    & \bm{p}^{\mathrm{min}} \le \bm{\alpha} \le \bm{\beta} \le \bm{p}^{\mathrm{max}}, \label{con2:prob1}
\end{align}
where the objective function in Eq.~\eqref{eq:obj_prob1} is an estimate obtained through cross-validation and cannot be expressed in closed form. 

For this reason, we consider applying the Nelder--Mead simplex method~\cite{nelder1965simplex,kochenderfer2019algorithms} to the black-box optimization problem (Eq.~\eqref{eq:obj_prob1}--\eqref{con2:prob1}). 
Note that this method can handle box constraints on the decision variables~\cite{ghiasi2008constrained} (e.g., $\bm{p}^{\mathrm{min}} \le \bm{\alpha} \le \bm{p}^{\mathrm{max}}$ and $\bm{p}^{\mathrm{min}} \le \bm{\beta} \le \bm{p}^{\mathrm{max}}$). 
Therefore, after incorporating the remaining constraints (i.e., Eq.~\eqref{con1:prob1} and $\bm{\alpha} \le \bm{\beta}$) into the objective function as penalty terms, we solve the following optimization problem:
\begin{align}
\underset{\bm{\alpha},\bm{\beta} \in \mathbb{R}^m}{\mathrm{maximize}} \quad 
& \mathrm{CV}_{\bm{X}}(\bm{\alpha},\bm{\beta}) - \lambda_1 g\Bigl(
    \sum_{j=1}^m (\beta_j - \alpha_j) - \Gamma \Bigr) - \lambda_2 \sum_{j=1}^m g\Bigl(\alpha_j - \beta_j\Bigr) \label{eq:obj_prob2}\\
\text{subject to} 
\quad & \bm{p}^{\mathrm{min}} \le \bm{\alpha} \le \bm{p}^{\mathrm{max}}, \quad \bm{p}^{\mathrm{min}} \le \bm{\beta} \le \bm{p}^{\mathrm{max}}, \label{con1:prob2} 
\end{align}
where $\lambda_1, \lambda_2 \in \mathbb{R}_{+}$ are sufficiently large penalty parameters, and $g(x) \coloneqq (\max\{x, 0\})^2$ is the penalty function for constraint violations\cite{audet2018mesh}.

\section*{Experimental results and discussion}
In this section, we evaluate the effectiveness of our methods for price bounds estimation through numerical experiments on synthetic datasets. 

\subsection*{Synthetic datasets}


Following the prior studies on price optimization~\cite{ito2017optimization,ito2018unbiased,ikeda2023prescriptive}, we created synthetic price--demand datasets (i.e., $\bm{X} \coloneqq \{(\bm{p}_i,\bm{d}_i) \mid i \in [n]\}$) using linear regression models. 
Here, the intercept term was sampled from a uniform distribution: $\theta_{j0}^{\star} \sim \mathrm{U}(m, 3m)$ for $j \in [m]$, and the regression coefficients were also sampled from uniform distributions: $\theta^{\star}_{j\ell} \sim \mathrm{U}(-2m, -3m)$ if $j = \ell$, and $\theta^{\star}_{j\ell} \sim \mathrm{U}(0, 3)$ otherwise, for $(j,\ell) \in [m] \times [m]$. 
The price of each item $j \in [m]$ for each instance $i \in [n]$ was sampled from a normal distribution: $p_{ij} \sim \mathrm{N}(0.8, 0.1^{2})$.

We then generated the demand for each item $j \in [m]$ for each instance $i \in [n]$ from the following linear regression models:
\begin{align} \label{eq:syn}
    d_{ij} = \theta_{j0}^{\star} + \sum_{\ell=1}^m \theta^{\star}_{j\ell} p_{i\ell} + \varepsilon_i \quad (i \in [n],j \in [m]),
\end{align}
where the Gaussian noise $\varepsilon_i \in \mathbb{R}$ for each instance $i \in [n]$ was sampled from a normal distribution: $\varepsilon_i \sim \mathrm{N}(0, \sigma^2)$. 
Here the noise variance $\sigma^2 \in \mathbb{R}$ was set according to the noise level, defined as follows: 
\begin{equation}
\delta \coloneq \sqrt{\frac{\sigma^2}{\mathrm{E}[d^2]}} 
= \sqrt{\frac{n m \sigma^2}{\sum_{i=1}^n \sum_{j=1}^m (d_{ij})^2}}. \quad 
\label{eq:delta}
\end{equation}

\subsection*{Methods for comparison}
We compared the performance of the following three methods for determining price bounds:
\begin{description}[leftmargin=!,labelwidth=\widthof{\bf CVNM:}]
    \item[Quantile:] Price bounds set to the $q$-quantiles of price data~\cite{ikeda2024interpretable} ($q \in \{60\%, 65\%, \ldots, 100\%\}$);
    \item[Bootstrap:] Our bootstrap estimation (Algorithm \ref{alg:bs});
    \item[Cross-validation:] Our cross-validation-guided black-box optimization, which solves problem (Eq.~\eqref{eq:obj_prob2}–\eqref{con1:prob2}) based on Algorithm~\ref{alg:kfoldCV} using the Nelder--Mead simplex method~\cite{nelder1965simplex,kochenderfer2019algorithms}.
\end{description}
All experiments were conducted on a Windows 11 computer equipped with an AMD Ryzen9 5900X CPU (3.70GHz, 12 cores) and 128 GB of RAM. 
To solve the price optimization problem (Eq.~\eqref{eq:obj}--\eqref{eq:constr_price}) and implement the Nelder--Mead simplex method, we used the Python SciPy library. 
The minimum and maximum feasible prices were set to $p^{\mathrm{min}}_j = 0.5$ and $p^{\mathrm{max}}_j = 1.1$, respectively, for all $j \in [m]$. 
In the bootstrap method, the number of bootstrap samples was set to $N^{\textrm{bs}} = 100$, and the critical-value parameter $\kappa \in \mathbb{R}$ was set according to the $\{60\%, 65\%, \ldots, 95\%, 99\%, 100\%\}$ confidence intervals. 
In the cross-validation method, the penalty parameters were set to $\lambda_1 = \lambda_2 = 1$, and the price range parameter was set as $\Gamma \in \{0.25, 0.50, \ldots, 3.00\}$. 

\subsection*{Evaluation metrics}

Let $\bm p^{\star} \in \mathbb{R}^m$ be the optimal prices corresponding to the ground-truth coefficient matrix, and $\hat{\bm p}(\bm \alpha, \bm \beta) \in \mathbb{R}^m$ be the optimal prices corresponding to the estimated coefficient matrix under the price bounds constraint as follows:
\begin{align}
    \bm p^{\star}
    \;\in\;
    \underset{\bm p \in \mathbb{R}^m} {\operatorname{argmax}}
    \Bigl\{ f(\bm p,\bm\Theta^{\star})
        \;\Bigm|\;
        \bm p^{\min} \le \bm p \le \bm p^{\max}
    \Bigr\}, \quad
    \hat{\bm p}(\bm \alpha, \bm \beta)
    \;\in\;
    \underset{\bm p \in \mathbb{R}^m} {\operatorname{argmax}}
    \Bigl\{ f(\bm p,\hat{\bm\Theta})
        \;\Bigm|\;
        \bm \alpha \le \bm p \le \bm \beta
    \Bigr\}. 
\end{align}
Following the prior studies on price optimization~\cite{ito2017optimization,ito2018unbiased,ikeda2023prescriptive}, we evaluated the quality of the solutions obtained by price bounds, by calculating the ratio of the total revenue to the maximum value based on the ground-truth coefficient matrix: 
\begin{equation}\label{eq:relative_revenue}
    \textbf{Relative revenue} \coloneqq \frac{f(\hat{\bm p}(\bm \alpha, \bm \beta),\bm\Theta^{\star})}{f(\bm p^{\star},\bm\Theta^{\star})}~(\le 1).
\end{equation}
We also examined the average width of the price range for each item to evaluate how well each method narrowed the price range: 
\begin{equation}\label{eq:average_width}
\textbf{Average price range width} \coloneqq \frac{1}{m}\sum_{j=1}^{m}(\beta_j -\alpha_j).
\end{equation}
Obtained results were averaged over 100 independent runs, and the standard errors are shown as error bars in the figures. 

\begin{figure}[thb!]
\centering
    \begin{minipage}[t]{0.45\linewidth}
        \centering
        \includegraphics[keepaspectratio, width=0.99\linewidth]{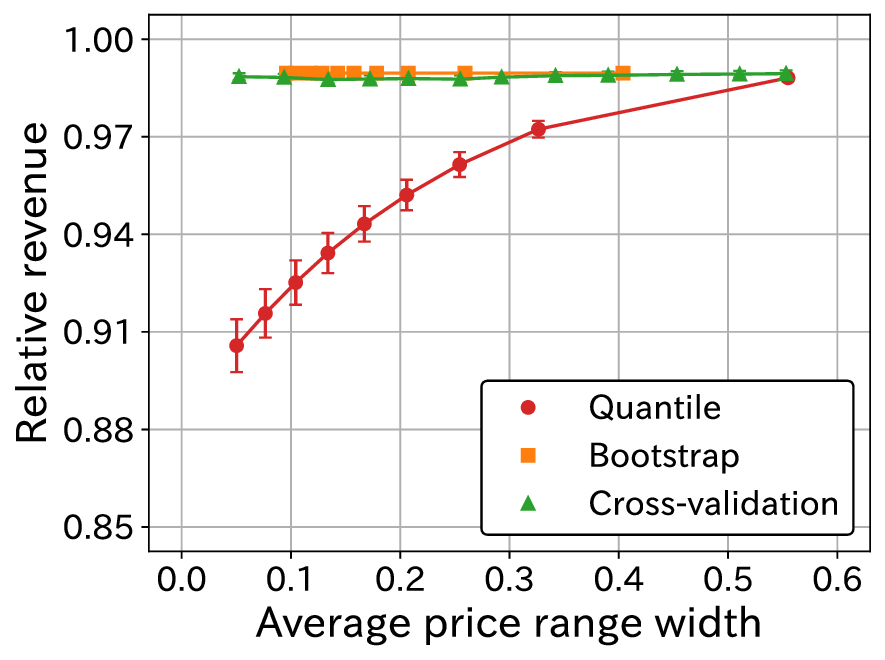}
        \subcaption{$\delta = 0.25$, $m=5$}
    \end{minipage}
    \begin{minipage}[t]{0.45\linewidth}
        \centering
        \includegraphics[keepaspectratio, width=0.99\linewidth]{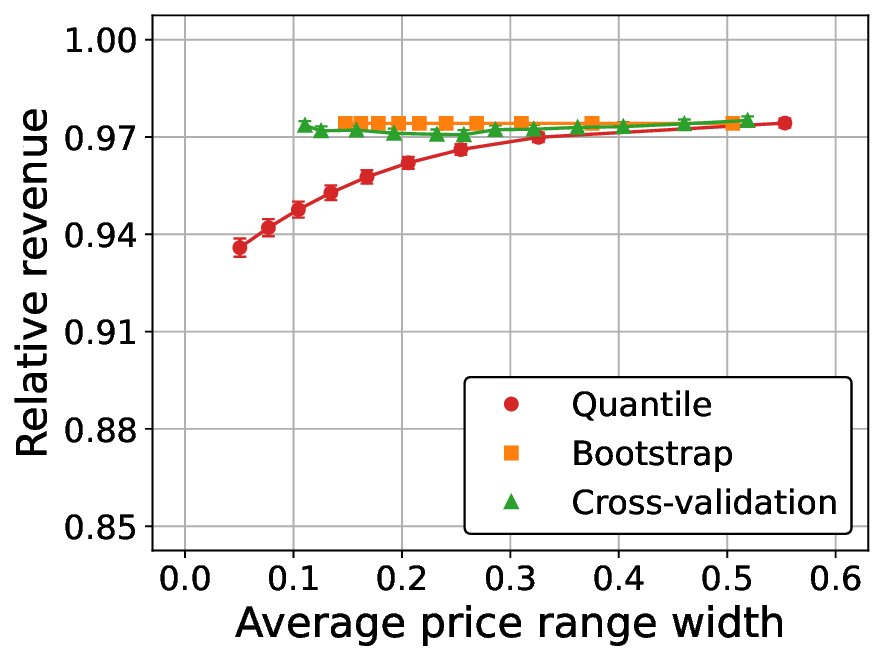}
        \subcaption{$\delta = 0.25$, $m=10$}
    \end{minipage}\\[5mm]
    \begin{minipage}[t]{0.45\linewidth}
        \centering
        \includegraphics[keepaspectratio, width=0.99\linewidth]{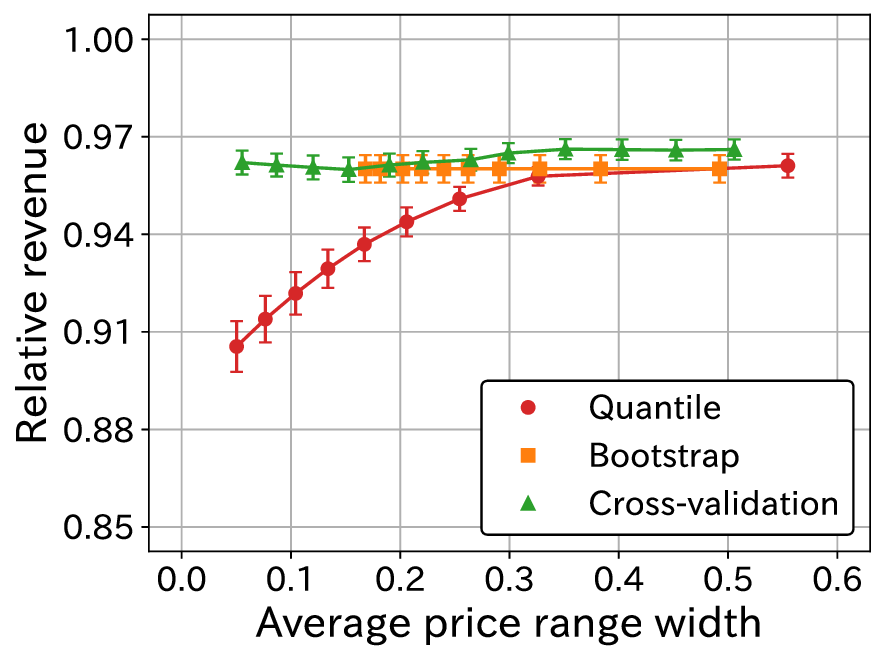}
        \subcaption{$\delta = 0.50$, $m=5$}
    \end{minipage}
    \begin{minipage}[t]{0.45\linewidth}
        \centering
        \includegraphics[keepaspectratio, width=0.99\linewidth]{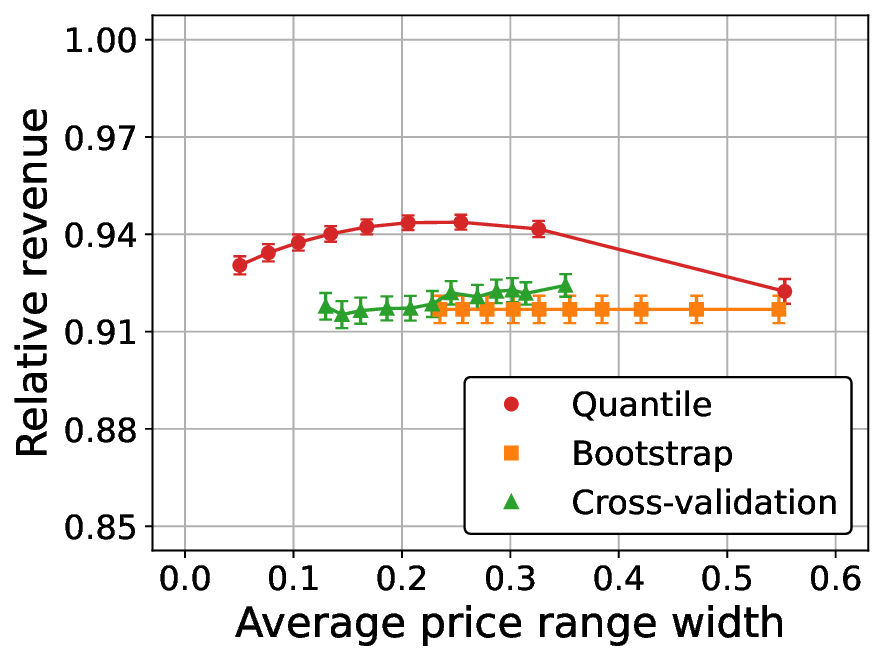}
        \subcaption{$\delta = 0.50$, $m=10$}
    \end{minipage} \\[5mm]
    \begin{minipage}[t]{0.45\linewidth}
        \centering
        \includegraphics[keepaspectratio, width=0.99\linewidth]{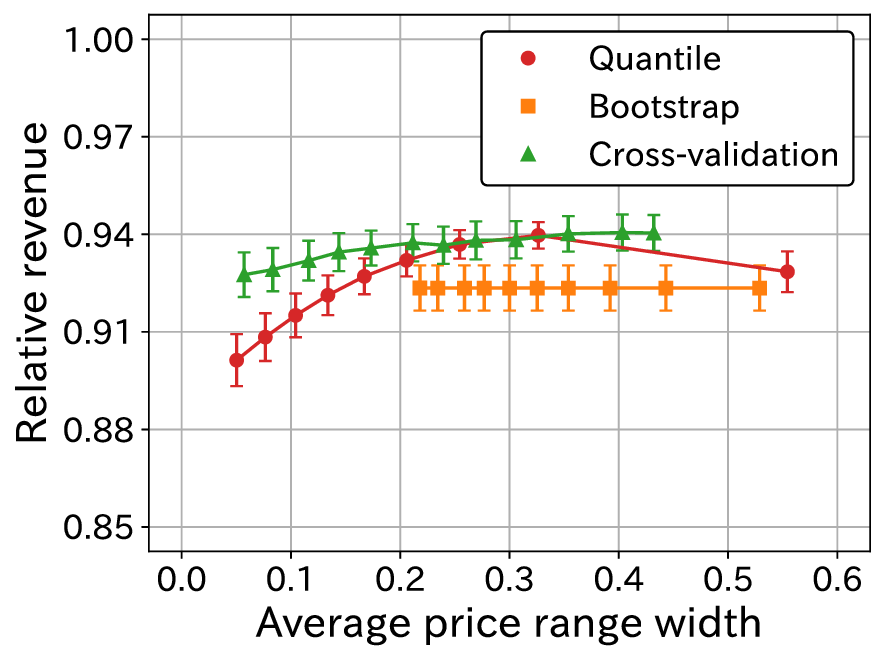}
        \subcaption{$\delta = 0.75$, $m=5$}
    \end{minipage}
    \begin{minipage}[t]{0.45\linewidth}
        \centering
        \includegraphics[keepaspectratio, width=0.99\linewidth]{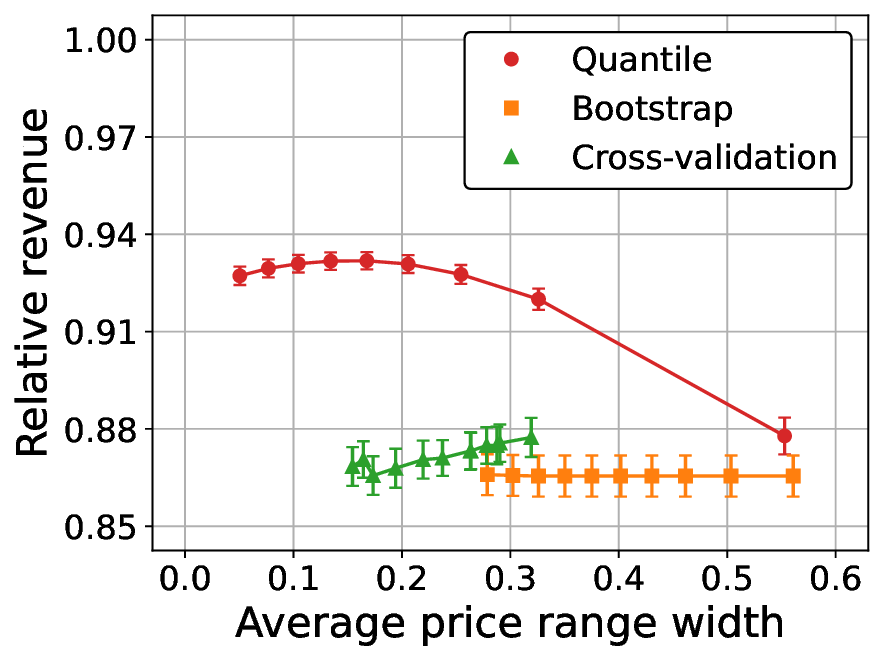}
        \subcaption{$\delta = 0.75$, $m=10$}
    \end{minipage}
\caption{Relative revenue as a function of the average price range width ($n = 300$)}
    \label{fig:noise_300}
\end{figure}

\begin{figure}[thb!]
\centering
    \begin{minipage}[t]{0.45\linewidth}
        \centering
        \includegraphics[keepaspectratio, width=0.99\linewidth]{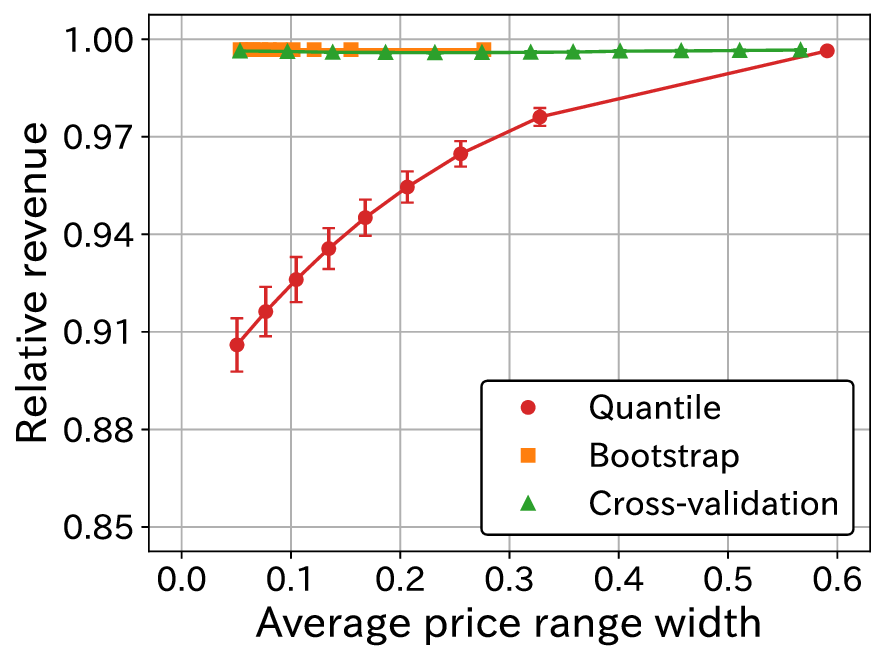}
        \subcaption{$\delta = 0.25$, $m=5$}
    \end{minipage}
    \begin{minipage}[t]{0.45\linewidth}
        \centering
        \includegraphics[keepaspectratio, width=0.99\linewidth]{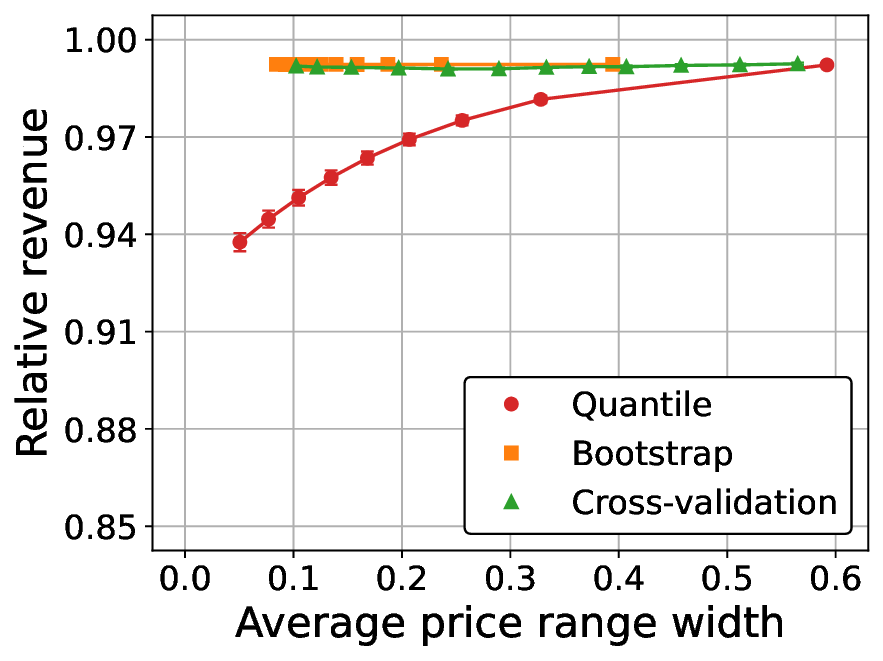}
        \subcaption{$\delta = 0.25$, $m=10$}
    \end{minipage}\\[5mm]
    \begin{minipage}[t]{0.45\linewidth}
        \centering
        \includegraphics[keepaspectratio, width=0.99\linewidth]{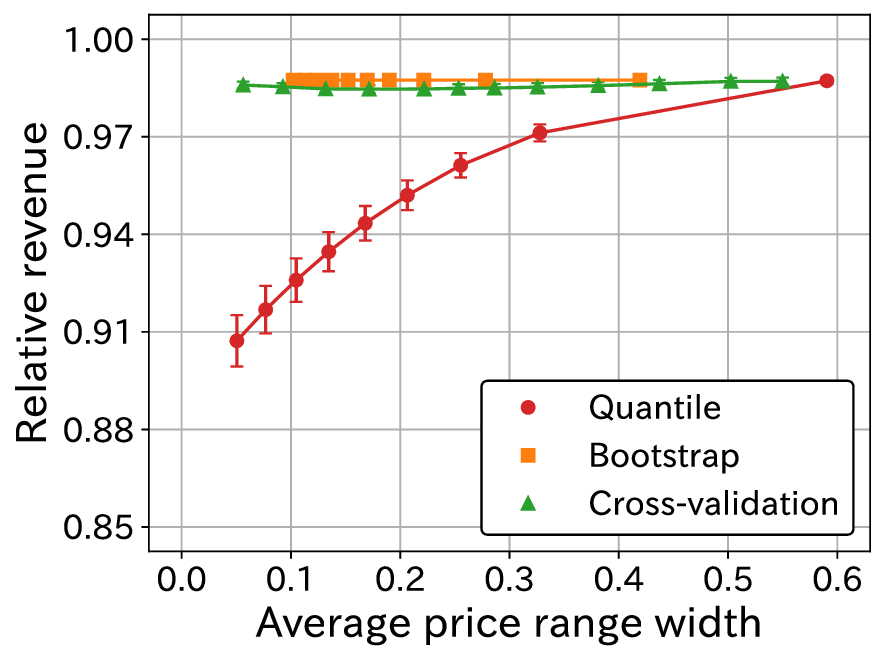}
        \subcaption{$\delta = 0.50$, $m=5$}
    \end{minipage}
    \begin{minipage}[t]{0.45\linewidth}
        \centering
        \includegraphics[keepaspectratio, width=0.99\linewidth]{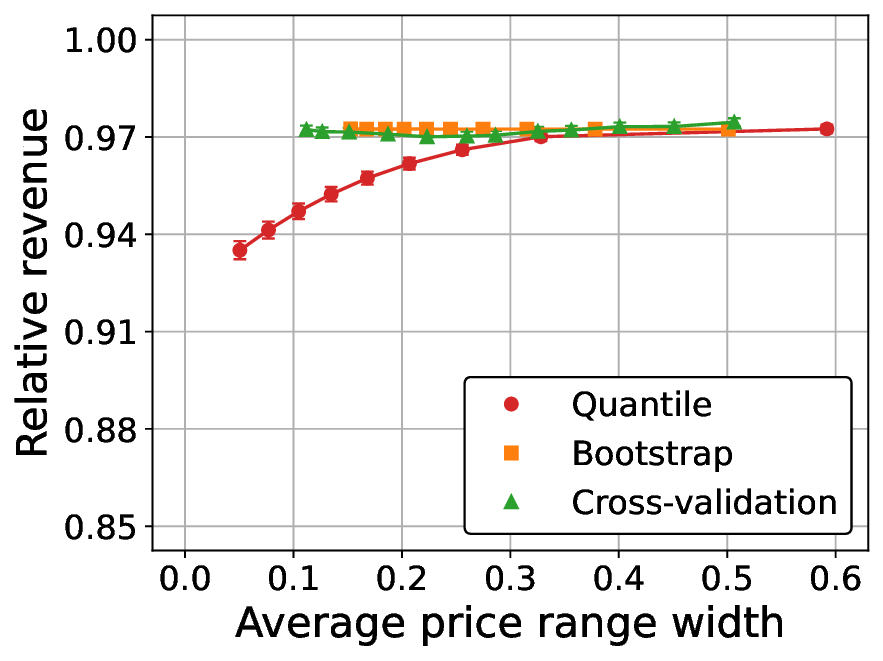}
        \subcaption{$\delta = 0.50$, $m=10$}
    \end{minipage}\\[5mm]
    \begin{minipage}[t]{0.45\linewidth}
        \centering
        \includegraphics[keepaspectratio, width=0.99\linewidth]{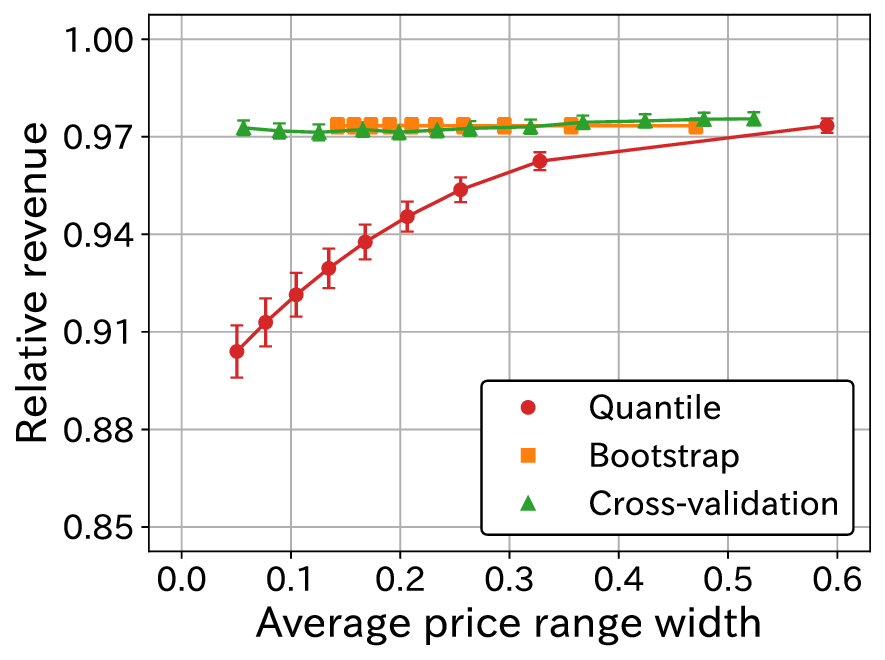}
        \subcaption{$\delta = 0.75$, $m=5$}
    \end{minipage}
    \begin{minipage}[t]{0.45\linewidth}
        \centering
        \includegraphics[keepaspectratio, width=0.99\linewidth]{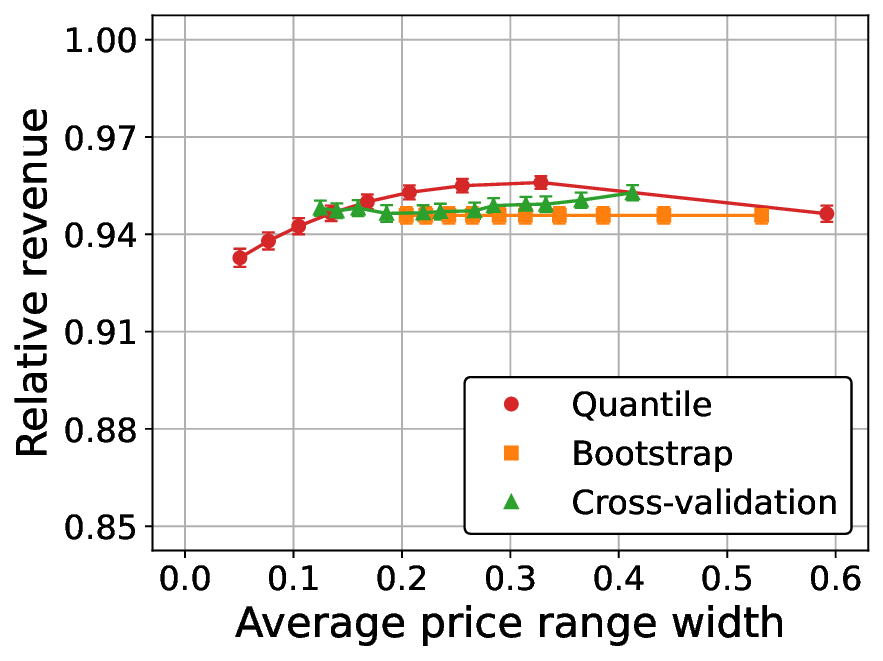}
        \subcaption{$\delta = 0.75$, $m=10$}
    \end{minipage}
\caption{Relative revenue as a function of the average price range width ($n = 1000)$}
    \label{fig:noise_1000}
\end{figure}

\subsection*{Results of total revenues}

Figures \ref{fig:noise_300} and \ref{fig:noise_1000} show the relative revenue (Eq.~\eqref{eq:relative_revenue}) as a function of the average price range width (Eq.~\eqref{eq:average_width}) for dataset sizes $n=300$ and $n=1000$, respectively. 
Here, we set $m \in \{5, 10\}$ as the number of items and $\delta \in \{0.25, 0.50, 0.75\}$ as the noise level.

We begin by focusing on the results for $n=300$ in Figure~\ref{fig:noise_300}. 
When $m=5$, our cross-validation method narrowed the price range while maintaining high revenues. 
Our bootstrap method also achieved high revenues similar to the cross-validation method for $\delta \in \{0.25, 0.50\}$, whereas its revenue deteriorated for $\delta = 0.75$. 
The conventional quantile method resulted in diminishing revenues as the price range narrowed.
When $m=10$, our bootstrap and cross-validation methods performed better than the quantile method for $\delta = 0.25$. 
For $\delta \in \{0.50, 0.75\}$, however, the quantile method improved revenues as the price range narrowed, outperforming our bootstrap and cross-validation methods.

We next move on to the results for $n=1000$ in Figure~\ref{fig:noise_1000}. 
When $m=5$, our bootstrap and cross-validation methods consistently outperformed the quantile method, successfully narrowing the price range without compromising revenues.
When $m=10$, our bootstrap and cross-validation methods also outperformed the quantile method for $\delta \in \{0.25, 0.50\}$, but were slightly inferior to the quantile method for $\delta = 0.75$.

These results indicate that our methods work well, especially for low-dimensional data with low noise.
The bootstrap method, which examines the distribution of optimal prices for bootstrap samples, and the cross-validation method, which maximizes ground-truth revenue estimates, allow us to limit the price range while maintaining high revenues.
On the other hand, the quantile method, which depends solely on historical price distributions, reduces revenues due to opportunity losses caused by excluding the optimal price from the price range. 
Moreover, as more data accumulates, the effectiveness of our methods can further improve compared to the quantile method.
In contrast, our methods may fail to perform well for high-dimensional noisy data. 
This suggests that simpler methods such as the quantile method may be more useful in such extremely difficult situations.

\begin{figure}[thb!]
\centering
    \begin{minipage}[t]{0.45\linewidth}
        \centering
        \includegraphics[keepaspectratio, width=0.99\linewidth]{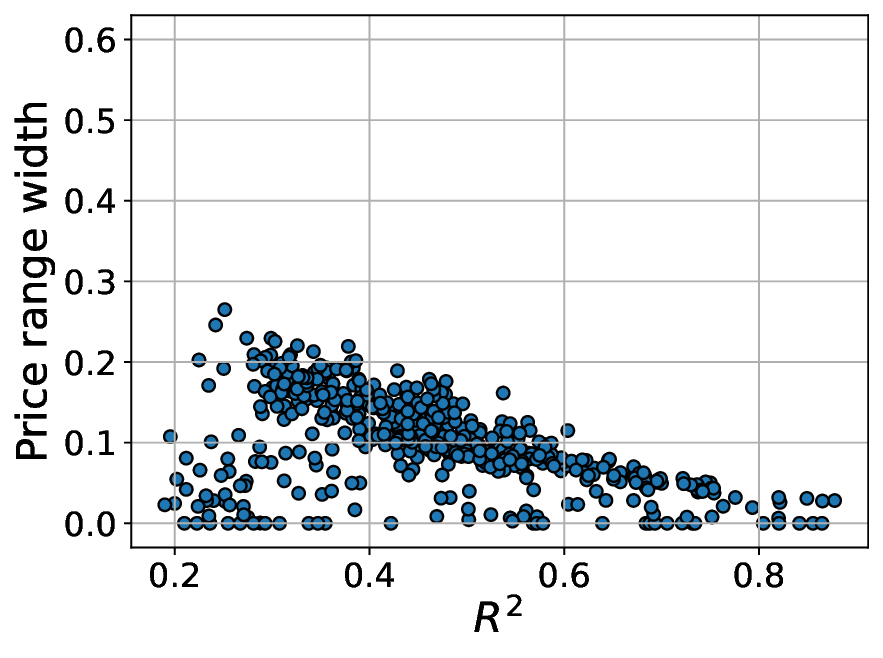}
        \subcaption{Bootstrap ($\delta=0.25$)}
    \end{minipage}
    \begin{minipage}[t]{0.45\linewidth}
        \centering
        \includegraphics[keepaspectratio, width=0.99\linewidth]{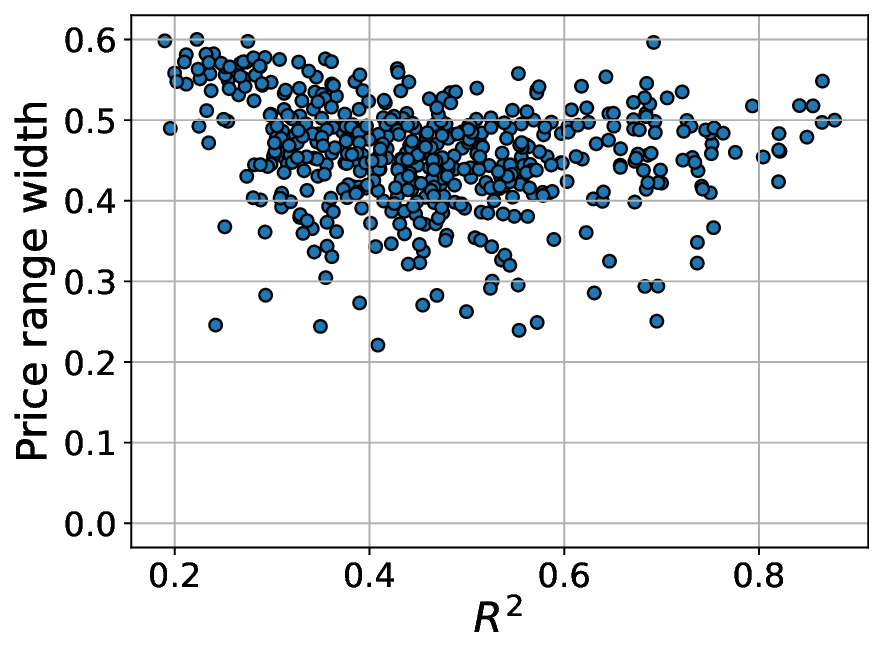}
        \subcaption{Cross-validation ($\delta=0.25$)}
    \end{minipage}\\[5mm]
    \begin{minipage}[t]{0.45\linewidth}
        \centering
        \includegraphics[keepaspectratio, width=0.99\linewidth]{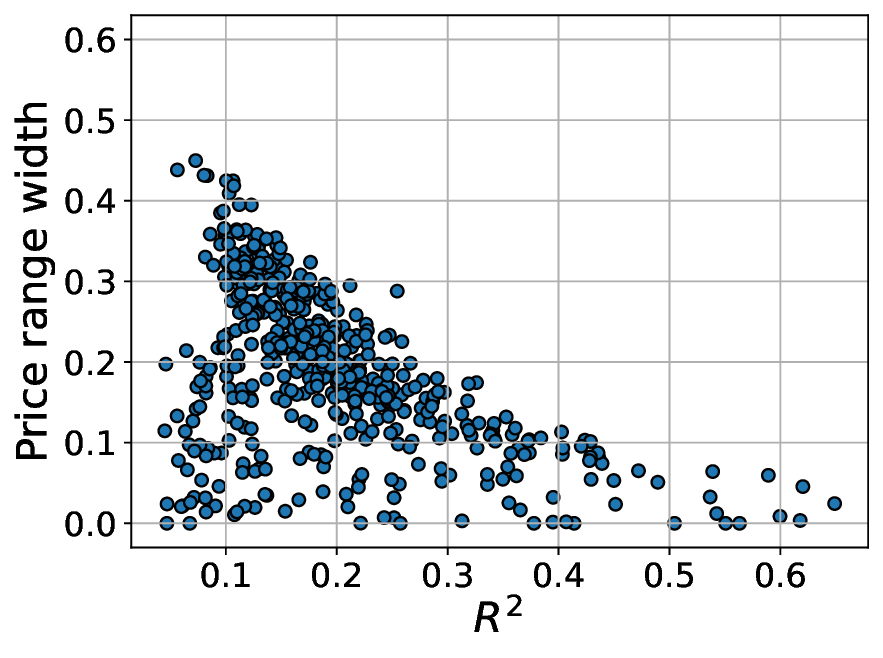}
        \subcaption{Bootstrap ($\delta=0.50$)}
    \end{minipage}
    \begin{minipage}[t]{0.45\linewidth}
        \centering
        \includegraphics[keepaspectratio, width=0.99\linewidth]{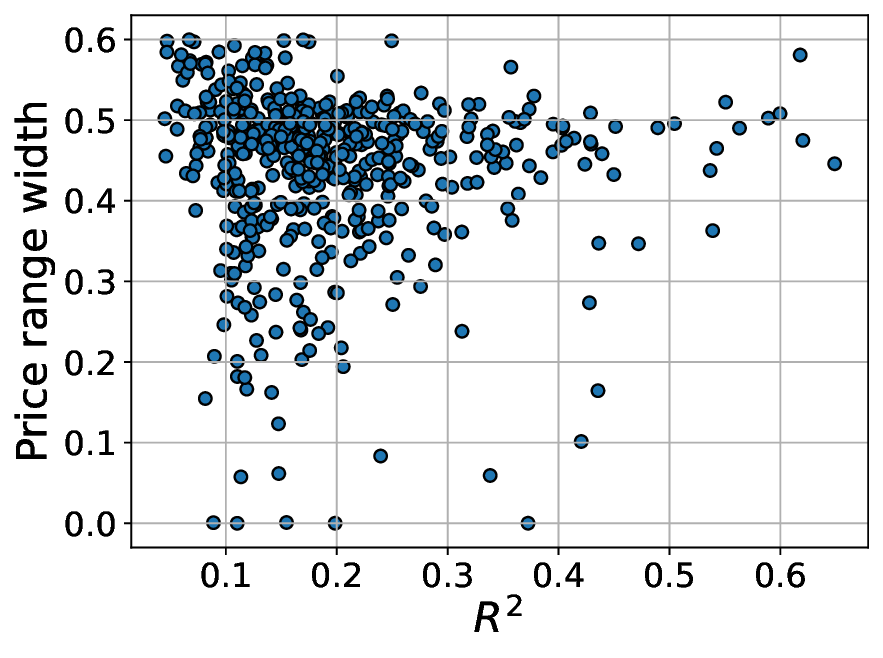}
        \subcaption{Cross-validation ($\delta=0.50$)}
    \end{minipage}\\[5mm]
    \begin{minipage}[t]{0.45\linewidth}
        \centering
        \includegraphics[keepaspectratio, width=0.99\linewidth]{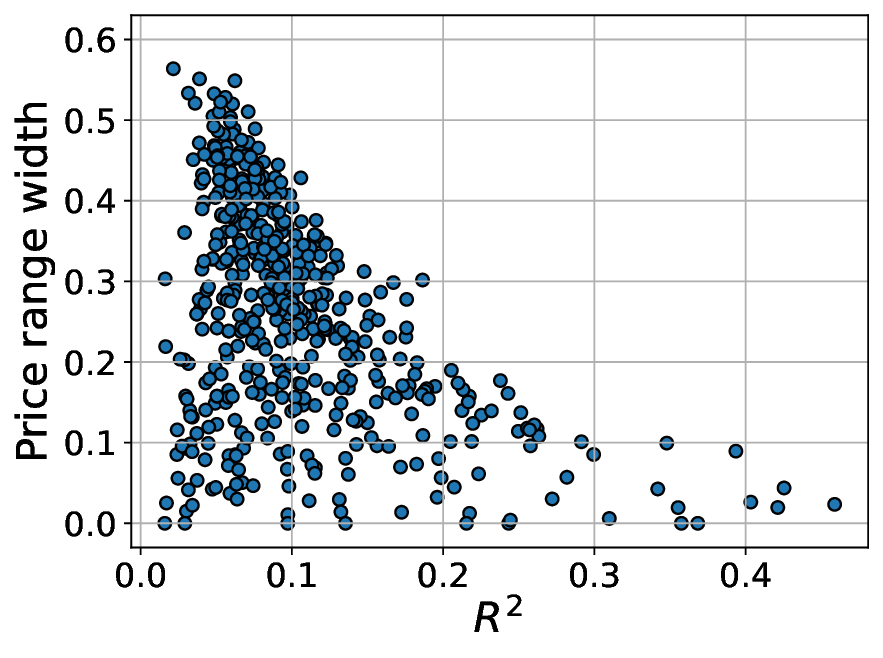}
        \subcaption{Bootstrap ($\delta=0.75$)}
    \end{minipage}
    \begin{minipage}[t]{0.45\linewidth}
        \centering
        \includegraphics[keepaspectratio, width=0.99\linewidth]{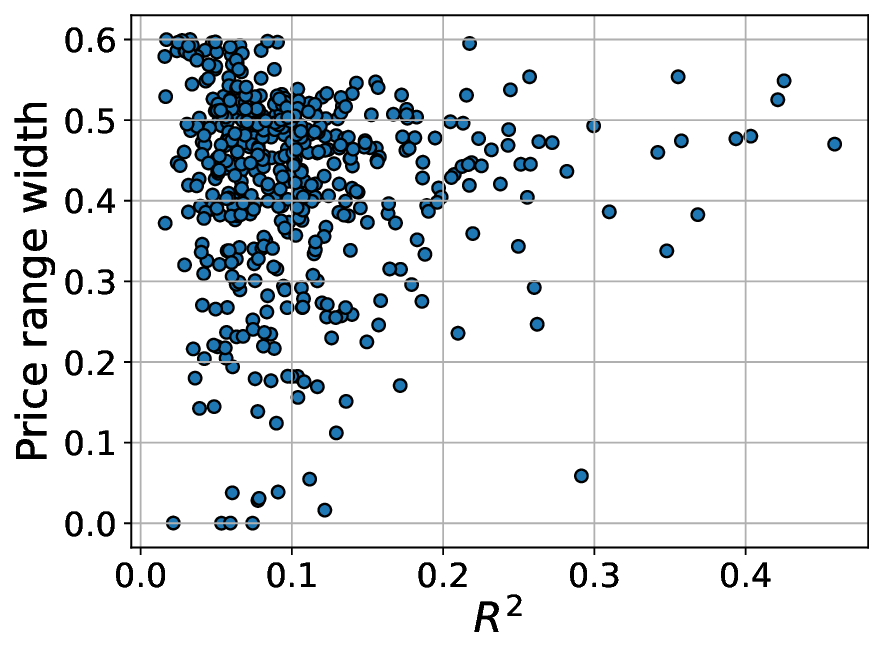}
        \subcaption{Cross-validation ($\delta=0.75$)}
    \end{minipage}
\caption{Scatter plot of demand forecasting accuracy and price range width for each item ($m=5$, $n=1000$, $\kappa=1.645$, $\Gamma=2.5$)}
    \label{fig:r2}
\end{figure}

\subsection*{Results of price range widths}

Figure \ref{fig:r2} shows scatter plots of the demand forecasting accuracy ($R^2$) and the price range width determined by our methods for each item.
Here, we set $m=5$ as the number of items, $n=1000$ as the dataset size, $\delta \in \{0.25, 0.50, 0.75\}$ as the noise level, $\kappa = 1.645$ as the bootstrap critical-value parameter, and $\Gamma = 2.5$ as the cross-validation price range parameter. 
Note that there are $5 \times 100 = 500$ points in each figure because 100 independent runs were conducted for the five-item price optimization problem.

The price ranges adjusted by the two methods showed contrasting trends depending on the accuracy of demand forecasts.
The bootstrap method set narrow price ranges overall for small $\delta$, and widened the price ranges as $\delta$ increased for items with low $R^2$.
When the demand forecasting accuracy is low, the optimal prices derived from bootstrap samples exhibit greater variance, resulting in wider price ranges.
In contrast, the cross-validation method set wide price ranges overall for small $\delta$, and narrowed the price ranges as $\delta$ increased for items with low $R^2$.
This is because the lower the accuracy of demand forecasts, the less reliable the optimal prices become, so narrowing the price range mitigates the risk of mispricing caused by unstable forecasts.

\begin{figure}[thb!]
\centering
    \begin{minipage}[t]{0.45\linewidth}
    \centering
    \includegraphics[keepaspectratio, width=0.99\linewidth]{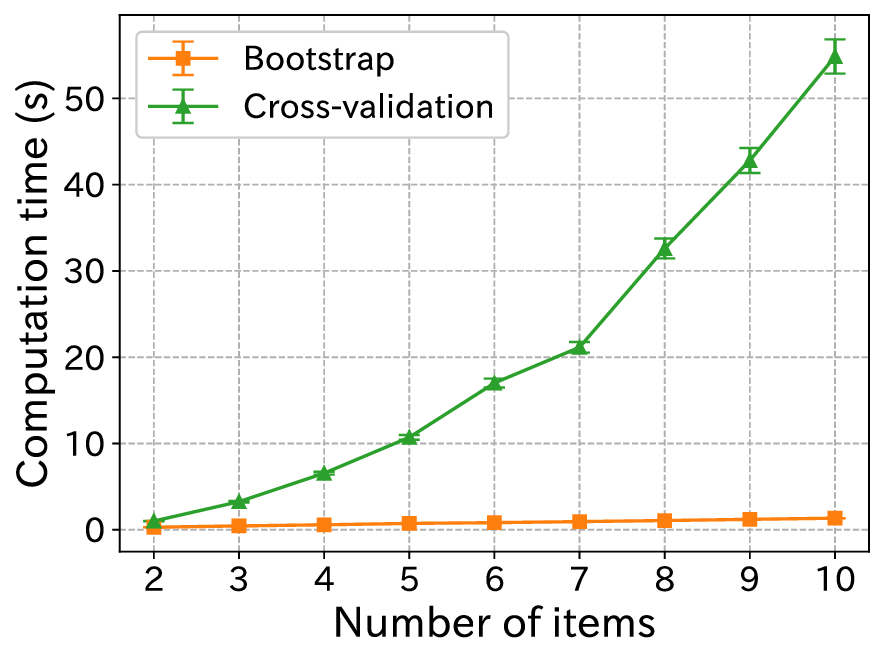}
    \subcaption{$n=300$}
    \end{minipage}
    \begin{minipage}[t]{0.45\linewidth}
    \centering
    \includegraphics[keepaspectratio, width=0.99\linewidth]{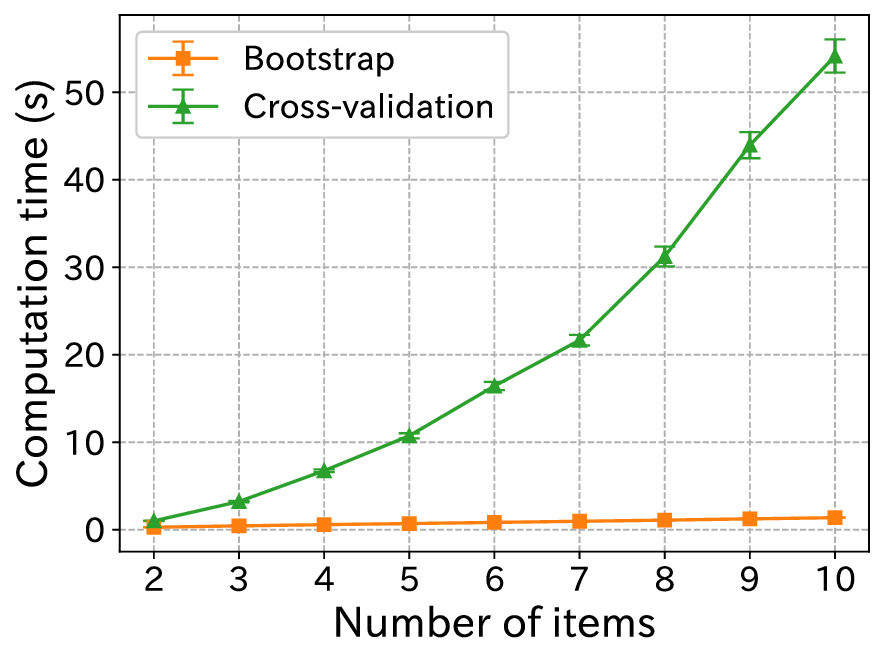}
    \subcaption{$n=1000$}
    \end{minipage}
\caption{Computation time in seconds as a function of the number of items ($\delta = 0.50$, $\kappa = 2.576$, $\Gamma = 3.00$)}
\label{fig:time}
\end{figure}

\subsection*{Results of computation times}

Figure \ref{fig:time} shows the computation time in seconds required by our methods as a function of the number of items $m \in \{2, 3, \ldots, 10\}$.
Here, we set $n \in \{300,1000\}$ as the dataset size, $\delta = 0.50$ as the noise level, $\kappa = 2.576$ as the bootstrap critical-value parameter, and $\Gamma = 3.00$ as the cross-validation price range parameter. 

Although the bootstrap method iterates demand forecasting and price optimization for each bootstrap sample, the computation time increased only slightly as the number of items increased.
On the other hand, when the number of items was large, the Nelder–Mead simplex method consumed a long time, increasing the overall computation time of the cross-validation method. 
Also note that both methods were largely unaffected by the dataset size.

\section*{Conclusion}
We proposed two methods to estimate price bounds for revenue maximization within the framework of prescriptive price optimization. 
The first method uses the bootstrap procedure~\cite{efron1992bootstrap,james2013introduction} to estimate confidence intervals for optimal prices. 
The second method uses the Nelder–Mead simplex
method~\cite{nelder1965simplex,kochenderfer2019algorithms} for black-box optimization to maximize total revenue estimated through $K$-fold cross-validation~\cite{ito2018unbiased}.

We conducted numerical experiments using synthetic datasets to validate the effectiveness of our methods for estimating profitable price bounds.
Our methods successfully narrowed down the price range while maintaining high revenues, especially for low-dimensional data with low noise.
As more data accumulated, the comparative advantage of our methods further increased.
Additionally, the price ranges adjusted by our two methods showed contrasting trends depending on the accuracy of demand forecasts as shown in Figure~\ref{fig:r2}.
The computation time for the bootstrap method increased only slightly with the number of items, whereas that for the cross-validation method was significantly affected by the item count.

A future direction of study will be to employ more advanced demand forecasting models in prescriptive price optimization.
While this research is based on linear regression models, it would be beneficial to incorporate predictive models that capture nonlinearity, such as regression trees~\cite{ferreira2016analytics,ikeda2023prescriptive}.
Another direction for future research will be to use the $K$-fold cross-validation procedure~\cite{ito2018unbiased} to tune various parameters involved in prescriptive price optimization. 
For example, potential applications include hyperparameter tuning of demand forecasting models and the construction of uncertainty sets in robust price optimization \cite{yabe2017robust,ikeda2025robust}.


\bibliography{ref}

\end{document}